# DISCUSSION OF "ANALYSIS OF VARIANCE—WHY IT IS MORE IMPORTANT THAN EVER" BY A. GELMAN

By Joop Hox and Herbert Hoijtink

*Utrecht University*

**Bayesian inference and fixed and random effects.** Professor Gelman writes "Bayesians see analysis of variance as an inflexible classical method." He adopts a hierarchical Bayesian framework to "identify ANOVA with the structuring of parameters into batches." In this framework he sidesteps "the overloaded terms fixed and random" and defines effects "as constant if they are identical for all groups in a population and varying if they are allowed to differ from group to group." Applying this approach to his first example (a Latin square with five treatments randomized to a $5 \times 5$ array of plots), variance components have to be estimated for row, column and treatment effects.

In our opinion, his approach provides an insightful connection between analysis of variance and hierarchical modeling. It renders an informative and easy to interpret display of variance components that is a nice alternative for traditional analysis of variance. However, we wonder whether sidestepping the terms fixed and random is always wise. Furthermore, currently his approach is rather descriptive, and does not contain truly Bayesian inference. Both points will be briefly discussed in the sequel.

To look into the question of fixed versus random and the use of hierarchical modeling, we carried out a small experiment. We constructed a dataset for the example in Section 2.2.2: 20 machines randomly divided into four treatment groups, with six outcome measures for each machine. We asked a statistician who is very skilled in multilevel analysis to analyze these data. The result was a hierarchical multivariate data structure with six outcomes nested within 20 machines, and the treatments coded as dummy variables at the machine level. Variance components were estimated for machines and measures. The treatment effects were tested by constraining all treatments to be equal and using a likelihood-ratio test.

Comparing this procedure with the discussion of this example in Gelman's paper shows that this is not what he had in mind. It certainly contradicts

---

the notion implied in Sections 3.2 and 3.3 that using hierarchical modeling, so to speak, automagically leads to a correct model. In fact, the multilevel analysis approach outlined above makes sense if we assume that the four treatments exhaust all treatments we are interested in. If we assume that there is a population of treatments, or that variations in implementation can lead to different outcomes, we can structure the data as a three-level model, with outcome measures nested within machines nested within treatments, and estimate a variance for the treatments. But even in this case one may ask if this variance is an interesting number to estimate. We would probably be more interested in the actual treatment effects, or in their differences.

Treating the treatment effects as fixed versus random requires knowledge about the actual design of the study, and a decision on how we should view these treatments. Our point here is that none of this comes automagically. We agree with Gelman that, once such decisions are made, the hierarchical modeling framework is both elegant and powerful. By way of illustration: all models discussed by Gelman for these data can be analyzed using the software MLwiN [Goldstein et al. (1998)]. Given the small sample size, maximum likelihood estimation is not attractive, but MLwiN includes a fully Bayesian inference option. So, at least one widely available multilevel program can be used to analyze these data correctly—after we have specified what we regard as "correct."

Our second point is that, to us, truly Bayesian inference is inseparably connected to the use of informative prior information (excluding the "I know nothing" kind of prior information) and the use of Bayes theorem to quantify the support in the data for different sources of prior information or competing models, that is, to compute posterior probabilities. Consider, for example, the five treatments in Gelman's first "25 plots Latin square" example. Even without the treatment labels, we can come up with several prior expectations that could be interesting in this context. Researcher $A$ evaluated the five treatments and came up with the following model for the treatment effects $\beta$:

$$M_A: \quad \beta_1 < \beta_2 < \beta_3 < \beta_4 < \beta_5.$$

Researcher $B$ has a different evaluation and renders the following prior expectation:

$$M_B: \quad \{\beta_1, \beta_4\} < \{\beta_2, \beta_3, \beta_5\}.$$

These models imply that the treatment effect is not (to use Gelman's terminology) a varying effect. Stated otherwise, the assumption that the five treatment effects come from the same distribution does not hold. It is also clear that the treatment effects are not constant, that is, equal for all treatment groups.

This problem can be solved by reinstating the term fixed effect and defining it as a varying effect with components that do not come from the same distribution. However, then the data cannot be analyzed in the framework proposed by Gelman. This does not imply that we disqualify his approach; we only want to stress again that there are situations in which the terms fixed and random effects (varying effects that do not and do come from the same distribution) are still appropriate.

Evaluation of models $A$ and $B$ (for simplicity ignoring the row and column effects of the Latin square design) is possible in a Bayesian framework. First of all, a prior distribution has to be specified for each model. If $\sigma^2$ denotes the residual variance of a one-way ANOVA with five groups, this prior could have the form

$$g(\theta|M_m) = g(\beta_1, \beta_2, \beta_3, \beta_4, \beta_5, \sigma^2|M_m) \propto \prod_{i=1}^{5} N(\beta_i|0, 1000)\chi^{-2}(\sigma^2|1, 10) I_{M_m},$$

where $i$ denotes the treatment, the indicator function has the value 1 if the treatment effects are in accordance with the restrictions imposed by model $m = A, B$ and 0 otherwise, and $N(\cdot)$ and $\chi^{-2}(\cdot)$ denote uninformative normal and scaled inverse chi-square distributions, respectively. Note that the resulting prior distributions are informative because the prior expectations formulated by researchers $A$ and $B$ are included using the indicator function. Note also that otherwise the prior is uninformative, and does not differ between treatment effect parameters and competing models. Subsequently, Bayes theorem can be used to compute the posterior probability of models $A$ and $B$:

$$P(M_m|y) \propto P(y|M_m)P(M_m),$$

where

$$P(y|M_m) = \int_\theta P(y|\theta)g(\theta|M_m)\,d\theta,$$

and $y$ is a vector containing the treatment effects for each of the 25 plots.

In our opinion this approach is truly Bayesian because prior knowledge is formalized in prior distributions and subsequently evaluated using posterior probabilities. This is lacking in Gelman's approach. His prior distribution is uninformative, and there are neither competing models nor different sources of prior information that are evaluated using posterior probabilities.

The remaining question is whether it is a problem that Gelman changes from fixed/random to constant/varying, and, whether it is a problem that his prior distribution is uninformative and that there is no inference in the sense that posterior probabilities are computed. Both the answers no and yes are possible. No, because the approach proposed is valuable in itself. Yes, because (as is hopefully illustrated by the examples given) fixed effects are not necessarily treated optimally if Professor Gelman's approach is used. Also yes, because the mainly descriptive framework presented by Professor Gelman can potentially be modified such that competing models/prior information can be evaluated in a Bayesian manner. Consider once more the Latin square example with row, column and treatment effects. A first model could state that the variance of the row and column effects is zero; a second model that the variance of the row and column effects is smaller than the variance of the treatment effects; and, a third model that the variance of the treatment effects is zero. Potentially, the Bayesian approach can be used to compute posterior probabilities for each of these models. The main problem is the specification of prior distributions. As has been illustrated, this is fairly easy for inequality constrained models. The construction of priors for the comparison of models with zero and nonzero variance components is less straightforward. Bluntly fixing a variance at zero for one model, and giving it an uninformative prior distribution for another model, will lead to an analogue of Lindley's paradox [Lee (1997), pages 130 and 131]. The solution here might be prior distributions based on training data [Berger and Pericchi (1996)], or informative prior distributions.

P.O. BOX 80140
NL-3508 TC UTRECHT
THE NETHERLANDS
E-MAIL: J.J.C.M.Hox@fss.uu.nl